\documentclass[a4paper,11pt]{article}
\usepackage{amsmath}
\usepackage{amssymb}  
\usepackage[utf8x]{inputenc}
\newtheorem{prop}{Proposition}

\newtheorem{lem}[prop]{Lemma}

\newtheorem{thm}[prop]{Theorem}
\newtheorem{cor}[prop]{Corollary}

\newtheorem{rem}[prop]{Remark}

\usepackage[T1]{fontenc}

\newcommand{\N}{\mathbb{N}}
\newcommand{\R}{\mathbb{R}}
\newcommand{\Z}{\mathbb{Z}}
\newcommand{\C}{\mathbb{C}}

\newcommand*\x{{\bf x}}

\title{\large{\textbf{SIMPLY AND TANGENTIALLY HOMOTOPY EQUIVALENT BUT NON-HOMEOMORPHIC HOMOGENEOUS MANIFOLDS}}}
\author{\small{SADEEB OTTENBURGER}}

\begin{document}

\maketitle

\begin{abstract}
For each odd integer $r$ greater than one and not divisible by three we give explicit examples of infinite families of simply and tangentially homotopy equivalent but
 pairwise non-homeomorphic closed homogeneous spaces with fundamental group isomorphic to $\Z/r$. 
As an application we construct the first examples of manifolds which possess infinitely many metrics of nonnegative 
sectional curvature with pairwise non-homeomorphic homogeneous souls of codimension three with trivial normal bundle, such that
their curvatures and the diameters of the souls are uniformly bounded.
These manifolds are the first examples of manifolds fulfilling such geometric conditions and they serve as solutions to a problem posed by I. Belegradek,
S. Kwasik and R. Schultz.
\end{abstract}

\section{Introduction}
In 1966 J. Milnor proved in \cite {M-66} that for each lens space $L$ with dim$(L)\geq 5$ and fundamental group of order five or greater 
than six there exist infinitely many pairwise distinct $h$-cobordisms $(W;L,L')$ over $L$ which are distinguished by the Whitehead torsion of 
$(W,L)$. In this case the Whitehead torsion was related to the Reidemeister torsion of the boundary components and Milnor's result implied the existence of 
infinitely many pairwise non-homeomorphic smooth closed manifolds in the $h$-cobordism class of $L$. These manifolds belong to the class of \textit{fake lens spaces} which are quotients of free actions of cyclic
groups on odd dimensional spheres.\\
C.T.C. Wall obtained classification results for the class of
fake lens spaces \cite {W-99} with fundamental group of odd order which led to infinite sequences of simply homotopy equivalent but pairwise non-homeomorphic 
smooth closed manifolds.\\
Let $M$ be a $(4k+3)$-dimensional oriented closed smooth manifold, where $k\geq 1$. S. Chang and S. Weinberger \cite{CW-03} implicitly proved that if $\pi_1(M)$ is not torsion-free
then there exist infinitely many closed smooth manifolds which are simply and tangentially homotopy equivalent to $M$ but pairwise non-homeomorphic.\\
The last two examples of infinitely many distinct manifolds sharing the same simple or tangential homotopy type were found by analyzing surgery exact sequences.
The existence of all these infinite series was proven implicitly.\\\\
The motivation of this work was to find \textit{explicit} infinite families of pairwise non-homeomorphic manifolds lying in the same simple and tangential homotopy type.\\\\
Let $\cal{L}$ be the set of total spaces $L^{p,q}$ of principal $S^1$-fibre bundles over $S^2\times S^2$ given by the first Chern class $px+qy$, where $x$ and $y$ 
are the standard generators of $H²(\cdot;\Z)$ of the first and the second factor of the base space respectively and $(p,q)\in\Z²$ not equal to $(0,0)$.
These manifolds form a subclass of the class of lens space bundles over $S^2$:\\\\
By $L^3_p$ we denote the standard 3-dimensional lens space with fundamental group isomorphic to $\Z/p$, equipped
with the induced metric from $S^3$.
The manifold $L^{p,q}$ is diffeomorphic to the total space of the $L^3 _p$-bundle over $S^2$ with clutching function
\[S^1\rightarrow \mathrm{Isom}(L^3_p),\ \ z\mapsto (x\mapsto z^qx).\]
\begin{thm}
Let $r,t$ be integers, where $r$ is odd, greater than one and not divisible by three. The set $\{L^{r,(t+kr)r}\vert k\in\Z\}$ consists of simply and tangentially homotopy equivalent but
pairwise non-homeomorphic manifolds on which $SU(2)\times SU(2)\times U(1)$ acts smoothly and transitively.
\end{thm}
In 1972 J. Cheeger and D. Gromoll \cite {CG-72} proved that any complete open Riemannian manifold with nonnegative sectional curvature is diffeomorphic 
to the total space of the normal bundle of a totally geodesic and totally convex submanifold, called a soul.\\
This fundamental structural result on nonnegatively curved complete Riemannian manifolds led to further questions and results concerning the existence of
 infinitely large families of codimension $(n-l)$-souls or the moduli space of complete metrics of nonnegative sectional curvature.\\
For example I. Belegradek \cite{B-03} constructed the first examples of manifolds admitting infinitely 
many nonnegatively curved metrics with mutually non-homeomorphic souls of codimension at least five. In V. Kapovitch, A. Petrunin and W. Tuschmann \cite{KPT-05}
found examples of manifolds admitting a sequence of nonnegatively 
curved metrics with pairwise non-ho\-meo\-mor\-phic souls of codimension at least eleven but with more geometric control, i.e. they gave uniform upper bounds for
the sectional curvature of the open manifolds and the diameter of the souls.
\\\\
J. Milnor proved in \cite{M-61} that homotopy equivalent 3-dimensional lens spaces share the property that taking the cartesian 
product with $\R^3$ yields diffeomorphic manifolds. But one cannot realize an infinite sequence of pairwise homotopy equivalent but non-homeomorphic
lens spaces because there are just finitely many lens spaces in each dimension for each finite cyclic group realized as the fundamental group.
\\
A motivation of \cite{BKS-09} was to construct infinitely many pairwise homotopy equivalent but non-homeomorphic
souls of the smallest possible codimension. They were able to construct infinitely many distinct codimension four souls but not of codimension less than four 
(Problem 4.8 (i) in \cite {BKS-09}).\\\\
As an application of Theorem 1 we obtain
\begin{thm}
Let $r,q$ be integers such that $r$ is odd, greater than one and not divisible by three. Then there exists a positive
constant $D$ independent of $r$ and $q$, such that
the manifold $L^{r,qr}\times \R^3=:M_{r,q}$ possesses an infinite sequence of metrics $\{g^i_{r,q}\}_i$ of nonnegative sectional curvature with pairwise non-homeomorphic
 souls $\{S^i_{r,qr}\}_i$, such that\[0\leq \mathrm{sec}(M_{r,q},g^i_{r,q} )\leq 1\textrm{ and }\mathrm{diam}(S^i_{r,q})\leq D,\]
where $SU(2)\times SU(2)\times U(1)$ acts smoothly, transitively and isometrically on these souls.\\
Moreover there is no manifold such that an infinite subset of manifolds in $\{S^i_{r,qr}\}_i$ can be realized as souls of codimension one or two with trivial normal bundle (Remark 21 (ii)).
\end{thm}
As a byproduct of the previous theorems there is 
\begin{prop}
There are infinitely many homotopy equivalent but pairwise non-diffeomorphic Riemannian non-simply 
connected 5-manifolds with $0\leq sec \leq 1$ and diameter $\leq D$, where $D$ is a positive constant.
\end{prop}
Proposition 3 yields counterexamples to a slightly relaxed version of a question by S.-T. Yau (see section 4).
\\\\
If we choose $q$ to be zero then $M_{r,q}=L^3_r\times S^2\times \R³$. Let $L$ be
a 3-dimensional lens space which is homotopy equivalent to $L^3_r$ then we know from \cite {M-61} that $L^3_r\times \R³$ and $L\times \R^3$ are diffeomorphic. Hence
 $L^3_r\times S^2\times \R³$ and $L\times S^2\times \R³$ are diffeomorphic and we obtain as a special case of Theorem 2.
\begin{cor} Let $r$ be as in the previous theorems and $L$ a lens space homotopy equivalent to $L_r^3$ then
 $L\times S^2\times \R³$ possesses an infinite sequence of metrics of nonnegative sectional curvature with pairwise non-homeomorphic
homogeneous souls such that the geometric implications in Theorem 2 hold.
\end{cor}
Let $N$ be a smooth manifold and $\mathfrak{R}^u_{sec\geq 0}(N)$ be the set of complete smooth metrics with topology induced by uniform smooth convergence and if we mod out
the action of $\mathrm{Diff}(N)$ on $\mathfrak{R}^u_{sec\geq 0}(N)$ via pullback we obtain the moduli space of complete smooth metrics on $N$ of nonnegative sectional curvature
$\mathfrak{M}^u_{sec\geq 0}(N)$. A soul is in general not unique but by V.A. Sharafutdinov \cite {Sh-79} any two souls of a metric can be mapped onto each other by a diffeomorphism of 
the total space. I. Belegradek, S. Kwasik and R. Schultz sharpened this result \cite [Theorem 1.4(iii)] {BKS-09} by proving that the function that assigns to a nonnegatively curved
complete metric on a manifold $N$ the diffeomorphism type of the pair $(N,S_g)$ is constant on connected components of $\mathfrak{M}^u_{sec\geq 0}(N)$. This gives another proof of the following known fact:
\begin{prop}
If $r$ is as in Theorem 1 then $\mathfrak{M}^u_{sec\geq 0}(L_r^3\times S^2\times \R^3)$ has infinitely many components.
\end{prop}
The proofs of Theorem 1 and Theorem 2 (see Lemma 19) imply
\begin{prop}Let $r$ be as in the previous theorems. There exist infinitely many pairwise inequivalent smooth transitive actions of
$SU(2)\times SU(2)\times U(1)\times \R^3$ on $L_r^3\times S^2\times \R^3$, i.e. they pairwisely don't differ by a self-diffeomorphism of $L_r^3\times S^2\times \R^3$,
whereas there exists only one smooth transitive action of $SU(2)\times SU(2)\times U(1)$ on $L_r^3\times S^2$ up to self-diffeomorphism.\\
\end{prop}
This work is structured as follows:\\\\
\textit{A homotopy classification:} 
Let $n,r$ be non-zero integers, where $r$ is chosen as in Theorem 1. 
We show that the manifolds in $\{L^{r,(t+kr)r}\}_k$ lie in the same simple and tangential homotopy type. Therefore we first give
a classification of the manifolds in $\cal{L}$ up to homotopy. 
We explain why the tangent bundle of each such manifold is stably parallelizable and why the Reidemeister torsion of them is trivial.
 These facts imply that the homotopy equivalences which may occur between these manifolds are tangential and simple. \\\\
\textit{Distinguishing homeomorphism types:} We introduce a diffeomorphism invariant, calculate it for the manifolds in question and conclude that the manifolds in
$\{L^{r,(t+kr)r}\}_k$ are pairwise non-diffeomorphic. In the smooth category of closed 5-manifolds non-diffeomorphic always implies non-homeomorphic which relies on 
the fact that there are no exotic smooth structures on $S^5$.
\\\\
\textit{Souls of codimension three:}
 We show that each manifold in $\cal {L}$ is diffeomorphic to the quotient of
$SU(2)\times SU(2)\times U(1)$ by an isometric 2-torus action. O'Neill's formula on Riemannian submersions ensures the existence of metrics of nonnegative sectional 
curvature on the corresponding manifolds.
\\\\
\textit{Remarks:} 
Here we explain why there is no manifold $M$ such that an infinite subset of manifolds in $\cal {L}$ can be realized as souls of $M$ with trivial normal bundle and of codimension
less than three.\\\\
\textbf{Acknowledgements}\\\\
I have to thank Diarmuid Crowley for helpful discussions. I am very grateful to Igor Belegradek 
for giving me useful advices on how to present the results. I want to thank
 Wilderich Tuschmann who drew my attention to this subject.
\section{A homotopy classification}
Let $L^{p,q}\in \cal {L}$. The Gysin sequence associated to 
the $S^1$-fibre bundle over $S^2\times S^2$ given by the first Chern class $px+qy$ implies that $\pi_1(L^{p,q})$ is isomorphic to $\Z/\textrm{gcd}(p,q)$.
Furthermore since $(p,q)\neq(0,0)$ we conclude that the universal covering space of $L^{p,q}$ is 
\[L^{\frac{p}{gcd (p,q)}\frac{q}{gcd (p,q)}}\]
and again the above mentioned Gysin sequence together with the Hurewicz theorem implies that $\pi_2(L^{p,q})$ is isomorphic to $\Z$.\\
Let $\Pi_{p,q}:L^{p,q}\rightarrow S^2\times S^2$ denote the projection map of the fibre bundle. The tangent bundle of $L^{p,q}$, $\tau_{L^{p,q}}$, 
is isomorphic to $\Pi_{p,q}^*(\tau_{S^2\times S^2})\oplus
\epsilon^1$, where $\epsilon^1$ denotes the trivial real line bundle over $L^{p,q}$. This description of $\tau_{L^{p,q}}$ implies via 
the pullback property of (stable) characteristic classes that its Stiefel-Whitney and Pontrjagin classes vanish.\\
From these facts we may deduce the following:
\begin{lem}
\begin{itemize}
 \item [(i)] The fundamental group of $L^{p,q}$ acts homotopically trivial on the universal covering space.
\item [(ii)] The tangent bundle of $L^{p,q}$ is stably trivial. Thus $L^{p,q}$ admits a spin structure and if the order
of $\pi_1(L^{p,q})$ is odd then the spin structure is unique.
 \item [(iii)] The universal covering space $\tilde{L}^{p,q}=L^{\frac{p}{gcd (p,q)}\frac{q}{gcd (p,q)}}$ is diffeomorphic to $S^2\times S^3$. Thus $\pi_2(L^{p,q})=\Z=\pi_3(L^{p,q})$.
 \end{itemize}
\end{lem}
\textbf{Proof.} (i) The deck-transformations sit inside a circle action. (ii) Since by definition $\Pi_{p,q}$ is onto and the tangent bundle of $S^2\times S^2$ is stably trivial we conclude that the classifying map of $\tau_{L^{p,q}}$,
 $c_{p,q}:L^{p,q}\rightarrow BO$, is zero homotopic and hence $\tau_{L^{p,q}}$ is stably parallelizable. (iii) Since the universal covering space of $L^{p,q}$ is spin it follows 
from a result of S. Smale \cite {S-62} that the diffeomorphism type only depends on $H_2(\tilde{L}^{p,q};\Z)$. The Gysin sequence of the $S^1$-fibre bundle
structure of $L^{\frac{p}{gcd (p,q)}\frac{q}{gcd (p,q)}}$ implies that $H_2(\tilde{L}^{p,q};\Z)=\Z$. Hence $\tilde{L}^{p,q}$ is
diffeomorphic to $S^2\times S^3$.
\hfill $\square$
\\\\
From now on let the manifolds in $\cal {L}$ be oriented by orienting the base and the fibre in the standard way and throughout this work if we speak about spin manifolds
we mean oriented manifolds given a certain spin structure.
\begin{thm}
 Let $r$ be as in Theorem 1 and $L^{p,q}$, $L^{p',q'}\in{\cal L}$ with $\pi_1(L^{p,q})\cong\pi_1(L^{p',q'})\cong \Z/r$.
 Furthermore let $(m,n)$, $(m',n')$ be pairs of integers such that $m\frac{q}{r}+n\frac{p}{r}=1=m'\frac{q'}{r}+n'\frac{p'}{r}$.
 Then $L^{p,q}$ and $L^{p',q'}$ are oriented homotopy equivalent if and only if there exist $s,s'\in(\Z/r)^*$, $\epsilon,\epsilon '\in \{\pm 1\}$ and $k,k'\in \Z/r$ such that
\[
s^3\frac{p,q}{r^2}\equiv  s'^3\frac{p',q'}{r^2}\mathrm{ mod \textrm{ }} r,\ \ \ \ \]
\[s (\epsilon m + k\frac{p}{r})(\epsilon n - k \frac{q}{r})\equiv  s' (\epsilon' m' + k'\frac{p'}{r})(\epsilon' n' - k' \frac{q'}{r})\mathrm{ mod \textrm{ }} r,\]
\[s^2(\frac{q}{r}(\epsilon m +k\frac{p}{r})-\frac{p}{r}(\epsilon n -k\frac{q}{r}))\equiv  s'^2(\frac{q'}{r}(\epsilon' m' +k'\frac{p'}{r})-\frac{p'}{r}
(\epsilon n' -k'\frac{q'}{r}))\mathrm{ mod \textrm{ }} r.
\]
\end{thm}
Before we give a proof of Theorem 8 we gather some further differential topological properties of the manifold $L^{p,q}$. 
\begin{lem}
The  \textit{Reidemeister torsion} in the sense of \cite [pp.404] {M-66} is defined for $L^{p,q}$ and it is trivial.
\end{lem}
\textbf{Proof.} This torsion invariant is defined for manifolds with finite cyclic fundamental group and the property that it acts trivially on 
the rational cohomology ring of the universal covering space which by Lemma 7 (i) and the fact that $\pi_1(L^{p,q})\cong \Z/\textrm{gcd}(p,q)$ is the case
 for $L^{p,q}$. Since $L^{p,q}$ is the 
total space of an $S^1$-fibre bundle over a 1-connected space it follows from \cite[Thm. B] {HKR-07} that its Reidemeister torsion is trivial.\hfill $\square$
\begin{lem}The second level of the Postnikov tower of $L^{p,q}$ is
\[L^\infty_{gcd(p,q)}\times \C P^\infty=:B_{gcd(p,q)}\]
 up to fibrewise homotopy equivalence, where $L_r^\infty$ is the 
 infinite dimensional lens space which is a $K(\Z/r,1)$ and $\C P^\infty$ is the infinite dimensional complex projective space which is a $K(\Z,2)$.
\end{lem}
\textbf{Proof.} As $\pi_1(L^{p,q})\cong \Z/\textrm{gcd}(p,q)$ the first level of the Postnikov decomposition of $L^{p,q}$ equals $L^\infty_{gcd(p,q)}$. 
Since $\pi_1(L^{p,q})$ acts trivially on the higher homotopy groups Postnikov theory implies that the second level of the Postnikov tower of $L^{p,q}$ is a
$\C P^\infty$-fibration over $L^\infty_{gcd(p,q)}$, which is a pullback of the principal $K(\pi_2(L^{p,q})\cong\Z, 2)$-fibration over $K(\pi_3(L^{p,q})\cong\Z,3)$.
By obstruction theory
the homotopy class of the classifying map of this fibration may be identified with an element of $H^3(L^\infty_{gcd(p,q)};\Z)$. But
this group is trivial.\hfill$\square$
\\\\
Let $r$ be as in Theorem 1. The next step is to relate a bordism group to homotopy classification:\\\\
The bordism group of our interest is $\Omega^{Spin}_5(B_r)$ which
is defined to be the set 
\[\left\{(M,f)\vert \begin{array}{cc}
	M\textrm{ a closed smooth} 
                            \textrm{ 5-dimensional spin manifold},\\
			      f:M\rightarrow B_r\end{array}\right \} \]
modulo an equivalence relation
which is given as follows: $(M,f)\sim(N,g)$ if there exists a 6-dimensional smooth spin manifold with boundary equal to the disjoint union of $M$ and $N$ and a map
 $F:W\rightarrow B_r$
which restricted to the boundary components is the map $f$, $g$ respectively.\\\\
The Postnikov decomposition of $L^{p,q}$ yields maps $f:L^{p,q}\rightarrow B_r$ which are 3-equivalences, i.e. they induce isomorphisms on the first 
and second homotopy groups. We call such a map a \textit{normal 2-smoothing}.\\
Let $f_0:L^{p,q}\rightarrow B_r$, $f_1:L^{p',q'}\rightarrow B_r$ be normal 2-smoothings
 then $(L^{p,q},f_0)$ and $(L^{p',q'},f_1)$ define elements in $\Omega^{Spin}_5(B_r)$. Assume
that $(L^{p,q},f_0)$ and $(L^{p',q'},f_1)$ represent the same element in $\Omega^{Spin}_5(B_r)$ then there exists a bordism $(W,F)$ between them as
described above.
\\\\
Let $L^{s,\tau}_6(\Z/r,S)$ be the set which consists of stable
equivalence classes of weakly based non-singular $(-1)$-quadratic forms over $\Z[\Z/r]=:\Lambda$, where stabilization goes by 
taking orthogonal sum with $(-1)$-hyperbolic forms. Weakly based means that there is a choice of an equivalence class of bases of the underlying free $\Lambda$-module,
where two bases are equivalent if the change of basis matrix has trivial Whitehead torsion. The $S$ between the brackets stands for the choice of a so called form parameter \cite {Ba-81}.
One can easily show that $L^{s,\tau}_6(\Z/r,S)$ is a group with group structure given by taking orthogonal sum. The simple Wall group $L^s_6(\Z/r,S)$ and $L^{s,\tau}_6(\Z/r,S)$ 
are related to each other by the following exact sequence:
\[0\rightarrow L^s_6(\Z/r,S) \stackrel{i}{\rightarrow} L^{s,\tau}_6(\Z/r,S)\stackrel{\tau}{\rightarrow}\textrm{Wh }(\Z/r),\]
where $i$ is just the canonical inclusion and the map $\tau$ sends stable equivalence classes of weakly based non-singular $(-1)$-quadratic forms 
over $\Lambda$ to the Whitehead torsion of the matrix representation of the adjoint of this quadratic form with respect to to the chosen weak equivalence class of basis.\\
One can associate to $(W,F)$ an element in $L^{s,\tau}_6(\Z/r,S(N_0\times I))$, where $S(N_0\times I)$ is explained below:
\\\\
Let us recall Wall's definition of a quadratic form on an even dimensional compact manifold. We equip $W$ with a base point and orient it at this point.
 Wall defines a skew-hermitian form $\lambda$ on the group of regular homotopy classes of immersions of $3$-dimensional spheres into $W$ which roughly speaking
 is given by transversal double point intersections which along two branches are joined with the base point such that this form takes values in $\Lambda.$ 
This form is called the \textit{equivariant intersection form} associated to $W.$ Similarly Wall assigns to each immersion $u$ an 
element $\mu (u)\in\frac{\Lambda}{\left\langle a+\bar a\right\rangle}$ which is given by self-intersection. If we compose $\mu$ with the quotient map onto 
$\frac{\Lambda}{\left\langle a+\bar a,1\right\rangle}$ we call the result $\tilde{\mu}(u).$\\\\
By \cite [Prop. 4] {K-99} we may assume that $F:W\rightarrow B_r$ is a 3-equivalence and we identify $\pi_1(W)$ with $\Z/r$.
Since $(W,N_i)$ is $2$-connected $\pi_3(W,N_i)$ and $H_3(W,N_i;\Lambda)$ are isomorphic under the relative Hurewicz homomorphism. Poincar\'e duality implies
 that $H_3(W,N_i;\Lambda)$ is the only possibly non-vanishing homology group of the pair $(W,N_i)$ and by \cite[Lemma 2.3.]{W-99} it follows that $H_3(W,N_i;\Lambda)$
 is a stably free $\Lambda$-module with a preferred equivalence class of $s$-basis (for a definition of $s$-basis see \cite [p.369]{M-66}). If we take connected sum 
of $W$ with $\#_k(S^3\times S^3)$ for some $k\in \N$ big enough we may assume that $H_3(W,N_i;\Lambda)$ is a free $\Lambda$-module (see \cite [p.723]{K-99}).
The intersection form $\lambda: H_3(W,N_0;\Lambda)\times H_3(W,N_1;\Lambda)\rightarrow \Lambda$ is unimodular which follows from Poincar\'e duality and \cite [Theorem 2.1.]{W-99} 
tells us that this form is even simple if $H_3(W,N_0;\Lambda)$ and $H_3(W,N_1;\Lambda)$ are equipped with preferred bases.
Let $K\pi_3(W)$ be $\mathrm{Ker}(F_*:\pi_3(W)\rightarrow \pi_3(B_r)),$ $K\pi_3(N_i)$ be $\mathrm{Ker}(f_{i*}:\pi_3(N_i)\rightarrow \pi_3(B_r))$ and $\mathrm{Im}K\pi_3(N_i)$ be
the image of $K\pi_3(N_i)$ under the homomorphism which is induced by the inclusion $N_i\hookrightarrow W.$
We claim that $\mathrm{Im}K\pi_3(N_0)= \mathrm{Im}K\pi_3(N_1):$\\
Assume there is an element $x\in \mathrm{Im}K\pi_3(N_0)$ that doesn't lie in $\mathrm{Im}K\pi_3(N_1).$ 
Then by the homotopy exact sequence associated to $(W,N_1)$ $x$ represents a non-trivial element in $\pi_3(W,N_1).$ As $\lambda$ is non-degenerate and $\pi_3(W)\rightarrow \pi_3 (W,N_0)$ is surjective there exists a $y\in\pi_3(W)$ such that $\lambda(x,y)\neq 0.$ But since $x$ is trivial in $\pi_3(W,N_0)$ we have $\lambda(x,y)=0$ which is a contradiction. By interchanging the roles of $N_0$ and $N_1$ the claim follows.\\
The map $\pi_3(W)\rightarrow \pi_3(W,N_i)$ induces an isomorphism
$\frac{K\pi_3(W)}{\mathrm{Im}K\pi_3(N_i)}\stackrel{\approx}{\rightarrow} \pi_3(W,N_i)$
which is seen with the help of the following diagram,\\
\[
\begin{array}[pos]{ccccccc}
&&\pi_{4}(B,W)\ \ \ \ &&&\\
&&\downarrow\ \ \ \ \searrow&&\\
&&\pi_3(N_i)\rightarrow \pi_3(W)\rightarrow\pi_3(W,N_i)\rightarrow 0  \\	
&&\searrow\ \ \ \ \ \ \ \downarrow\ \ \ \ \ \ \ \ \ \ \ \ \ \ \ \ \ \\
&& \pi_3(B)\ \ \ \ \ \ \ \ \ \ \\
&&\downarrow\ \ \ \ \ \ \ \\
&&0,\ \ \ \ \ \ \ 
\end{array}
\]
where the diagonal maps are surjective.
\\\\
As promised above we explain what $S$ stands for: By $S(W)$ we denote the subgroup of $\Lambda$ which projects onto the image of $\mu$ restricted to $\mathrm{Im}K\pi_3(N_0)$.
Since $\mathrm{Im}K\pi_m(N_0)=\mathrm{Im}K\pi_m(N_1)$ $S(W)$ equals $S(N_0\times I)=S(N_1\times I)$ and we define $S$ to be $S(W)\oplus \Z$.
\\
If we equip $\frac{K\pi_3(W)}{\mathrm{Im}K\pi_3(N_0)}$ with the basis which is induced by the preferred basis on $H_3(W,N_0;\Lambda)$
then by the  map $\pi_3(W)\rightarrow \pi_3(W,N_i)$ which comes from the inclusion $W\hookrightarrow (W,N_i)$, $\lambda$ induces 
the form \[\bar{\lambda}:\frac{K\pi_3(W)}{\mathrm{Im}K\pi_3(N_0)}\times \frac{K\pi_3(W)}{\mathrm{Im}K\pi_3(N_0)}\rightarrow \Lambda.\]
This is a unimodular skew-hermitian form and 
$(\bar{\lambda},\tilde{\mu})$ represents an element $[(\bar{\lambda},\tilde{\mu})]=:\Theta(W,F)$ in $L^{s,\tau}_6(\Z/r,S)$ which indeed doesn't
depend on the choice of a normal bordism $(W,F)$.
\begin{rem}
 In \cite {Ot-11} an analysis of $\Theta(W,F)$ led to a classification of the corresponding manifolds up to diffeomorphism.
\end{rem}
For every $\Theta\in L_{6}^{s,\tau}(\Z/r,S)$ there exists a normal 2-smoothing $(L',f')$ and a bordism $(W',F')$ between $(L^{p,q},f)$
 and $(L',f')$ such that $\Theta (W',F')=\Theta.$ Moreover $(L',f')$ is up to $s$-cobordism completely determined by $\Theta.$ This is proved in 1
\cite [Theorem 5.8] {W-99} for the case of simple normal $n$-smoothings into a finite Poincar\'e complex. The proof
there extends literally to our situation. Furthermore $L^{p,q}$ and $\Theta\cdot L^{p,q}$ are homotopy equivalent.\\
These observations have the following consequence:\\
If two normal 2-smoothings $(L^{p,q},f)$ and $(L^{p',q'},g)$ represent the same element in $\Omega^{Spin}_5(B_r)$ then they are automatically homotopy equivalent.\\\\
Now we show that the converse is also true for spin 5-manifolds with finite cyclic fundamental group of odd order not divisible by 3 and with vanishing first Pontrjagin class.
\begin{prop}
Let $r$ be as in Theorem 1 and $N,M$ be closed smooth oriented spin 5-manifolds with vanishing first Pontrjagin classes and $f,g$ maps
 from $N,M$ to $B_r$ respectively Then $(N,f)$ and $(M,g)$ represent the same element in $\Omega^{Spin}_5(B_r)$ if and only if
\[f_*[N]=g_*[M]\in H_5(B_r;\Z).\]
\end{prop}
Since the condition in the previous proposition is a purely homotopical one we have the following
\begin{cor}
 Let $L$, $L'\in\cal {L}$ and $r$ as in Theorem 1 such that $\pi_1(L)\cong\Z/r\cong\pi_1(L')$ then $L$, $L'$ are homotopy equivalent if and only if
there exist normal 2-smoothings $f:L\rightarrow B_r$, 
$g:L'\rightarrow B_r$ such that $(L,f)$ and $(L',g)$ represent the same element in $\Omega^{Spin}_5(B_r)$.
\end{cor}
\textbf{Proof of Proposition 9.}
The entry $E_{a,b}^2$ of the $E_2$-term of the \textit{Atiyah Hirzebruch Spectral Sequence} (AHSS) for computing
\[\Omega_\star^{Spin}(L_r^{\infty}\times \C P^{\infty}=B_r)\] is $H_a(L_r^\infty\times \C P^\infty;\Omega^{Spin}_b(pt.))$. 
Since $r$ is odd the $E^2$-term for $a+b\leq 6$ looks as follows:
\\\\
$\begin{array}{cccccccccc}
6& &0 & & & & & & &  \\
5& &0 &0 & & & & & & \\
4& &\Z &\Z/r &\Z & & & & &  \\
3& &0 &0 &0 &0 & &   \\
2& &\Z/2 &0 &\Z/2 &0 &\Z/2 &   \\
1& &\Z/2 &0 &\Z/2 &0 &\Z/2 & 0  \\
0& &\Z &\Z/r &\Z &\Z/r &\Z &(\Z/r)^{3} &\Z \\
 & & & & & & & & &  \\
 & &0 &1 &2 &3 &4 &5 &6 \ \ \  .
\end{array}$\\\\
From \cite[p.7]{T-93} we know that the differentials in $E^{2}_{a,b}(L_r^\infty)$ respectively $E^{2}_{a,b}(\C P^{\infty})$ from the first to the 
second row are the dual of the Steenrod square $Sq^2:H^*(\cdot; \Z/2)\rightarrow H^{*+2}(\cdot; \Z/2)$ precomposed with the reduction
map in homology and the differentials from the second to the third row are the dual of $Sq^2$. 
The differentials $d_5:E_{6,0}^5(L_r^\infty)\rightarrow E_{1,4}^5(L_r^\infty)$ and
 $d_5:E_{6,0}^5(\C P^{\infty})\rightarrow E_{1,4}^5(\C P^{\infty})$ are trivial. This observation, the exterior product structure
 of $E_{s,t}^{r} (L_r^\infty\times \C P^\infty)$ induced by $E^{r}_{a,b}(L_r^\infty)$ and $E^{r}_{c,d}(\C P^\infty)$ and the fact that 
the differentials obey the Leibniz rule imply that
 $d_5:E_{6,0}^5(L_r^\infty\times \C P^\infty)\rightarrow E_{1,4}^5(L_r^\infty\times \C P^\infty)$ is trivial. Thus $E^{\infty}_{a,b}(L_r^\infty\times \C P^\infty)$ 
equals $E^{5}_{a,b}(L_r^\infty\times \C P^\infty)$ for $a+b=1,3,5$:\\\\
$\begin{array}{ccccccccc}
5& &0 & & & &  \\
4& & &\Z/r & &&  \\
3& &0 & &0 & &  \\
2& &&0 & &0 & \\
1& &0 &&0 & &0 \\
0& &&\Z/r & &(\Z/r)^{2} &&(\Z/r)^{3}& \\
 & & & & && & \\
 & &0 &1 &2 &3 &4 &5  \ \ \ .
\end{array}$\\\\
Hence
\begin{eqnarray*}h_1:\Omega _1^{Spin}(L_r^{\infty}\times \C P^{\infty})&\rightarrow& H_1(L_r^{\infty}\times \C P^{\infty};\Z),\\    \left[(S,l)\right ]&\mapsto& l_*[S]\end{eqnarray*}
is an isomorphism.\\
Let $K$ be a Kummer surface equipped with its standard orientation. We know from \cite {M-63} that $K$ generates $\Omega_4^{Spin}(pt.)$.
The construction of the AHSS and its $\infty$-term imply the following extension problem:
\begin{eqnarray}\Omega _1^{Spin}(L_r^{\infty}\times \C P^{\infty})\stackrel{\mu_K}{\hookrightarrow} 
\Omega _5^{Spin}(L_r^{\infty}\times \C P^{\infty})\stackrel{h_5}{\rightarrow} H_5(L_r^{\infty}\times \C P^{\infty};\Z)\rightarrow 0,\end{eqnarray}
where 
$\mu_K([(S,l)])= \textrm{$[(K\times S,l\circ \textrm{pr}_{S})]$}$ and $h_5([(N,g)])= g_*[N].$
Let $S^1$ be equipped with the standard orientation and $i:S^1\rightarrow L_r^{\infty}$ be the inclusion of $S^1$ as the 1-skeleton of $L_r^{\infty}$. 
The fact that $h_1$ is an isomorphism implies that $(S^1,i)$ represents a generator of $\Omega _1^{Spin}(L_r^{\infty}\times \C P^{\infty})$.\\
Let $v_1$ be a generator of $H^1(L_r^{\infty}\times \C P^{\infty};\Z/r).$ We claim that precomposing the homomorphism
\begin{eqnarray*}n_1:\Omega _5^{Spin}(L_r^{\infty}\times \C P^{\infty})&\rightarrow&\Z/r,\\
\left[(N,g)\right]&\mapsto&\left\langle \rho_r(p_1(N))g^*(v_1),[N]_{\Z/r}\right\rangle\end{eqnarray*}
with $\mu_K$ is an isomorphism:
\\The Kuenneth theorem implies that $(n_1\circ\mu_K)([S^1,i])$ equals
\[\left\langle \rho_r(p_1(K)),[K]_{\Z/r}\right\rangle\left\langle i^*(v_1),[S^1]_{\Z/r}\right\rangle,\]
where $\rho_r$ is the mod-$r$-reduction homomorphism in cohomology.
It is clear that $\left\langle i^*(v_1),[S^1]_{\Z/r}\right\rangle$ is a generator of $\Z/r$.
By the Hirzebruch signature theorem it is known that $\left\langle \frac{p_1}{3}(K),[K]\right\rangle$ equals $\textrm{sign}(K)$
 which is $-16.$ Thus \[\left\langle \rho_r (p_1(K)),[K]_{\Z/r}\right\rangle \equiv -48 \textrm{ mod }  r\]
is a generator of $\Z/r$ if and only if $\mathrm{gcd}(r,48)=1$.
But $\mathrm{gcd}(r,48)=1$ if and only if $r$ is odd and not divisible by 3 which is the case by assumption.\\\\
Thus $(n_1\circ\mu_K)([S^1,i])\in(\Z/r)^*$ and we can compose $n_1$ with an appropriate isomorphism $\alpha$ from $\Z/r$ to 
$\Omega _5^{Spin}(L_r^{\infty}\times \C P^{\infty})$ such that $\alpha\circ n_1$ is a splitting of (1). We conclude that
 $\Omega_5^{Spin}(L_r^{\infty}\times \C P^{\infty})\cong (\Z/r)^{4}$ and $(N,f)$, $(M,g)$ 
represent the same element in $\Omega_5^{Spin}(L_r^{\infty}\times \C P^{\infty})$ if and only if $n_1(N,f)=n_1(M,g)\textrm{ and }f_*[N]=g_*[M].$
Since the first Pontrjagin classes of the underlying manifolds are trivial the proposition follows.\hfill $\square$
\\\\
In order to prove Theorem 8 we need the following
\begin{lem} Let $r$ be as in Theorem 1 and $N$, $N'$ be smooth closed spin 5-manifolds with vanishing first Pontrjagin classes and $\pi_1(N)\cong\pi_1(N')\cong \Z/r$,
where $\pi_1$ acts trivially on $\pi_2$.
Then there exist normal 2-smoothings $(N,f)$ and $(N',f')$ representing the same element in $\Omega_5^{Spin}(B_r)$
if and only if there exist generators $v\in H^1(N;\Z/r)$, $v'\in H^1(N';\Z/r)$ and $z\in H^2(N;\Z)$, $z'\in H^2(N';\Z)$ which project to generators of
 $\frac{H^2(N;\Z)}{\mathrm{torsion}}$ and $\frac{H^2(N';\Z)}{\mathrm{torsion}}$ respectively such that
\begin{eqnarray*}
1)&&\left\langle v(\beta_r(v))^2,[N]_{\Z/r}\right\rangle\equiv\left\langle v'(\beta_r(v'))^2,[N']_{\Z/r}\right\rangle\mathrm{ mod \textrm{ }} r;\\
2)&&\left\langle v\beta_r(v)z_r,[N]_{\Z/r}\right\rangle\equiv\left\langle v'\beta_r(v')z'_r,[N']_{\Z/r}\right\rangle\mathrm{ mod \textrm{ }} r;\\
3)&&\left\langle vz_r^2,[N]_{\Z/r}\right\rangle\equiv\left\langle v'z_r'^2,[N']_{\Z/r}\right\rangle\mathrm{ mod \textrm{ }} r;
\end{eqnarray*}
where $z_r$ and $z'_r$ are the$\mathrm{mod}$-$r$-reductions of $z$ and $z'$ respectively.
\end{lem}
\textbf{Proof.}
Proposition 12 tells us that $(N,f)$ and $(N',f')$ represent the same element in $\Omega_5^{Spin}(L_r^{\infty}\times \C P^{\infty})$
if and only if $f_*[N]=f'_*[N']$. We observe that 
$H_5 (L_r ^\infty\times \C P^\infty; \Z)\cong H_5 (L_r ^\infty\times \C P^\infty; \Z /r)\cong H^5 (L_r ^\infty\times \C P^\infty; \Z /r).$
A basis of $H^5 (L_r ^\infty\times \C P^\infty; \Z /r)$ is given by $v_1 y_r^2,$ $v_1 (\beta _r (v_1)) y_r,$ $v_1 (\beta _r (v_1))^2,$
where $v_1$ is a generator of $H^1(L_r ^\infty\times \C P^\infty; \Z /r),$ $y_r$ is a generator of $H^2(L_r ^\infty\times \C P^\infty; \Z /r)$ which comes from the $\mathrm{mod}$-$r$-reduction of the standard generator $y$ of $H^2( \C P^\infty; \Z )$ and $\beta_r$ is the $\mathrm{mod}$-$r$ Bockstein homomorphism.
We see that $f_*[N]=f'_*[N']$ if and only if $f_*[N]_{\Z /r}=f'_*[N']_{\Z /r}.$ But this is the 
case if and only if
\begin{eqnarray*}
\left\langle f^*(v_1 (\beta _r (v_1))^2 ),[N]_{\Z /r}\right\rangle&=&\left\langle f'^*(v_1 (\beta _r (v_1))^2 ),[N']_{\Z /r}\right\rangle,\\														 
								\left\langle f^*(v_1 (\beta _r (v_1)) y_r),[N]_{\Z /r}\right\rangle&=&
\left\langle f'^*(v_1 (\beta _r (v_1)) y_r),[N']_{\Z /r}\right\rangle,\\																
\left\langle f^*(v_1 y_r^2),[N]_{\Z /r}\right\rangle&=&\left\langle f'^*(v_1 y_r^2),[N']_{\Z /r}\right\rangle.
\end{eqnarray*}
We finish the proof by replacing $f^*(v_1)$ by $v$, $f'^*(v_1)$ by $v'$, $f^*(y_r)$ by $z_r$ and $f'^*(y_r)$ by $z'_r.$
\hfill $\square$\\\\
Let $j:S^1\hookrightarrow L^{p,q}$ be the inclusion of the fibre which preserves the chosen orientation of the fibre. Furthermore let $m,n\in\Z$ such that 
$m\frac{q}{r}+n\frac{p}{r}=1.$ The Gysin sequence associated to $\Pi_{p,q}$ implies that $H^2(L^{p,q};\Z)\cong \Z/r\oplus \Z,$ where $\frac{p}{r}\Pi_{p,q}^*(x)+\frac{q}{r}\Pi_{p,q}^*(y)$ is 
a generator of the torsion part and $m\Pi_{p,q}^*(x)-n\Pi_{p,q}^*(y)$ is a generator of a $\Z $-summand. By $\alpha$ we denote the \textit{preferred generator} of $H^1(L^{p,q};\Z/r)$ which is 
characterized by the following property: \[\left\langle j^*(\alpha),[S^1]_{\Z/r}\right\rangle=1.\]
Let $v_1$ and $z$ be the standard generators of $H^1(L_r^{\infty};\Z/r)$ and $H^2(\C P^{\infty};\Z)$ respectively.\\
A normal 2-smoothing $f:L^{p,q}\rightarrow B_r$ is up to homotopy uniquely determined by 
\begin{eqnarray*}
f^*(v_1)&=&s\alpha,\\ f^*(z)&=&\epsilon(f)(m\Pi_{p,q}^*(x)-n\Pi_{p,q} ^* (y))+k(f,m,n)(\frac{p}{r}\Pi_{p,q}^*(x)+\frac{q}{r}\Pi_{p,q}^*(y)),
\end{eqnarray*}
where $s$ is a unit in $\Z /r$, $\epsilon(f)\in\{\pm 1\}$ and $k(f,m,n)\in \Z /r$.\\
We denote the set of triples $\{(s,\epsilon ,k)\vert s\in(\Z/r)^*, \epsilon \in\{\pm1\},k\in\Z/r\}$ by $T.$
\begin{lem}
Fixing a choice of $m,n\in\Z$ such that $m\frac{q}{r}+n\frac{p}{r}=1$ then there is a 1-1 correspondence between the set $H$ of homotopy classes of 2-smoothings
 of $L^{p,q}$ and $T,$ where the bijection is given as follows: 
\begin{eqnarray*}{\cal C} :&H\rightarrow& T, \\
  &[f]\mapsto& (s(f),\epsilon (f), k(f,m,n)).
\end{eqnarray*}
\end{lem}
\textbf{Proof.} It is clear that ${\cal C}$ is injective. We claim that ${\cal C}$ is also surjective, i.e. for fixed $m,n\in\Z$ as above any triple
$(\epsilon, s,k)\in T$ has a preimage under ${\cal C}$. We can write $f$ as $f_1 \times f_2$. The homotopy class $[f_1]$ of $f_1$ can be seen as an 
element in $H^1(L^{p,q};\Z /r)\cong \Z /r$.\\
Any self-automorphism of $\Z /r$ is given by a unit $s$ of $\Z /r$, $(1\mapsto s)$. Further there is a 1-1 correspondence between self-automorphisms of
 $\pi_1(L_r^{\infty})(\cong \Z /r)$ and homotopy classes of self-maps of $L_r^{\infty}$. Thus the homotopy classes of self-maps of $L_r^{\infty}$ correspond to
 self-automorphisms of $H_1(L_r^{\infty};\Z)$.\\
Now let $f$ be a self-map of $L_r^{\infty}$, then "naturality" of the Universal Coefficient Theorem (UCT) implies:
\begin{eqnarray*}
f^*:H^1(L_r^{\infty};\Z)\cong \mathrm{Hom }(H_1(L_r^{\infty};\Z),\Z /r)&\rightarrow &\mathrm{Hom }(H_1(L_r^{\infty};\Z),\Z /r),\\
h                                                            &\mapsto&f^*(h)=h\circ f_*.
\end{eqnarray*}
This means that $g^*:H^1(L_r^{\infty};\Z)\rightarrow H^1(L_r^{\infty};\Z)$ is an automorphism if and only if $g_* :H_1(L_r^{\infty};\Z)\rightarrow H_1(L_r^{\infty};\Z)$ is 
an automorphism. Thus the set of self-maps of $L_r^\infty$ that induce self-automorphisms on $\pi_1(L_r^\infty)$ is in 1-1 correspondence to $(\Z /r)^*$ which 
corresponds bijectively to self-automorphisms of $H^1(L_r^{\infty};\Z)$ and via Whitehead's theorem bijectively 
to the homotopy classes of self-homotopy equivalences of $L_r^{\infty}$. Hence by precomposing the $f_1$ in $f=f_1 \times f_2: L^{p,q}\rightarrow L_r^{\infty}\times \C P^\infty$ 
by a suitable self-homotopy equivalence one can realize any $s$ as in the above sense.\\
Without effecting $f_1$, we show that the homotopy class of $f_2:L^{p,q}\rightarrow \C P^{\infty}$, that realizes $(\epsilon, k)$, induces an isomorphism on $\pi_2$. Therefor
 we gather some facts:\\
a) From the Leray-Serre spectral sequence for the fibration $\widetilde{L}^{p,q}\stackrel{\mathrm{Pr}}{\rightarrow}L^{p,q}\rightarrow L_r^{\infty}$ we get

\[\begin{array}{ccccccccc}
0 & \longrightarrow & E_{20}^\infty & \longrightarrow & H^2(L^{p,q};\Z)& \stackrel{u}{\longrightarrow}& E_{02}^\infty & \longrightarrow & 0\\
          &  & \parallel               &            &    \parallel        & &                       \parallel         & & \\
0 & \longrightarrow  & H^2(L_r^\infty;\Z )& \longrightarrow & H^2(L^{p,q};\Z) & \stackrel{u}{\longrightarrow} & H^2(\widetilde{L}^{p,q};\Z ) & \longrightarrow & 0
\end{array} \]
with $u=\mathrm{Pr}^*$. Hence $\mathrm{Pr}^*:H^2(L^{p,q};\Z)\rightarrow H^2(\widetilde{L}^{p,q};\Z)$ is surjective with kernel isomorphic to $\Z /r$. \\\\
b) Further there exists the following commutative diagram:

\setlength{\unitlength}{1cm}
\begin{picture}(10,3)
\put(4,2.5){$\widetilde{L}^{p,q}$}
\put(4,0.5){$L^{p,q}$}
\put(6.5,0.5){$\C P^\infty$}

\put(4.2,2.3){\vector(0,-1){1.4}}
\put(4.8,0.6){\vector(1,0){1.6}}
\put(4.5,2.3){\vector(3,-2){2.2}}

\put(3.75,1.5){$\mathrm{Pr}$}
\put(5.9,1.5){$\tilde{f_2}$}
\put(5.3,0.2){$f_2$}
\end{picture}
\\
Applying the $\Z$-cohomology functor $H^2(\cdot ;\Z )$ we get the following commutative diagram:

\setlength{\unitlength}{1cm}
\begin{picture}(10,4)
\put(3,2.5){$H^2(\widetilde{L}^{p,q};\Z)$}
\put(3,0.5){$H^2(L^{p,q};\Z)$}
\put(7,0.5){$H^2(\C P^\infty;\Z)$.}

\put(4,0.9){\vector(0,1){1.4}}
\put(6.8,0.6){\vector(-1,0){1.8}}
\put(7.1,0.9){\vector(-3,2){2.2}}

\put(3.45,1.5){$\mathrm{Pr}^*$}
\put(6.2,1.7){$\tilde{f_2}^*$}
\put(5.6,0.2){$f^*_2$}
\end{picture}
\\
Since $H^2(\widetilde{L}^{p,q};\Z)\cong \Z,\ H^2(L^{p,q};\Z)\cong \Z\oplus \Z/r$ and $H^2(\C P^\infty;\Z)\cong \Z$ we see that $\tilde{f}_2^*$ must be an isomorphism.\\\\
c) Again "naturality" of the UCT implies that the set of homotopy classes of
maps from $\widetilde{L}^{p,q}$ to $\C P^\infty$ that induce isomorphisms on $H^2(\cdot; \Z )$ equals the set of homotopy classes of maps that induce 
isomorphism on $H_2 (\cdot;\Z)$.\\\\
d) Applying the Hurewicz theorem one sees that a map between simply connected CW-complexes that induces isomorphism on $H_2 (\cdot;\Z)$ also induces 
isomorphism on $\pi_2 (\cdot)$.\\\\
So by b) $\tilde{f}_2^*$ is an isomorphism on $H^2(\cdot ;\Z )$. Then c) implies that $\tilde{f}_2$ induces an isomorphism on $H_2 (\cdot;\Z)$
which via d) implies that $\tilde{f}_2$ induces isomorphism on $\pi_2 (\cdot)$ and thus by b) $f_2$ induces isomorphism on $\pi_2 (\cdot)$. Hence  ${\cal C}$ is surjective.
\hfill $\square$\\\\ 
\textbf{Proof of Theorem 8.} By Lemma 14 $(L^{p,q},f)$ and $(L^{p',q'},f')$ represent the same element in $\Omega_5^{Spin}(B_r)$
if and only if 
\begin{eqnarray}
\left\langle f^*(v_1)(\beta_r(f^*(v_1)))^2,[N]_{\Z/r}\right\rangle\equiv\left\langle f'^*(v_1)(\beta_r(f'^*(v_1)))^2,[N']_{\Z/r}\right\rangle\mathrm{mod\textrm{ }}r,\\
\left\langle f^*(v_1z_r)\beta_r(f^*(v_1)),[N]_{\Z/r}\right\rangle\equiv\left\langle f'^*(v_1z_r)\beta_r(f'^*(v_1)),[N']_{\Z/r}\right\rangle\mathrm{mod\textrm{ }}r,\\
\left\langle f^*(v_1z_r^2),[N]_{\Z/r}\right\rangle\equiv\left\langle f'^*(v_1z_r'^2),[N']_{\Z/r}\right\rangle\mathrm{mod\textrm{ }}r.
\end{eqnarray}
\textbf{Notation:} By $x_c, y_c\in H^2(S^2\times S^2; \Z/c)$ we denote the elements of $H^2(S^2\times S^2; \Z/c)$ which are the mod $c$ reductions
of the elements that come from the standard generator of $H^2(\cdot;\Z)$ of the first factor and the second factor of $S^2\times S^2$ respectively.\\\\
We know that $f^*(\beta_r (v_1))=\beta_r(f^*(v_1))=s\beta_r(\alpha).$
Hence in order to compute the Kronecker products above we have to understand what $\beta _r (\alpha)$ is in terms of $\Pi_{p,q} ^*(x)$ and $\Pi_{p,q} ^*(y)$, i.e. 
$\beta _r (\alpha)=b_1 \Pi_{p,q} ^*(x_r)+ b_2 \Pi_{p,q} ^*(y_r)$
for some $b_1,b_2\in \Z/r.$ The cohomological structure of $L_r^{\infty}\times \C P^{\infty}$ implies that $\beta _r (\alpha)$ lies in the image of $\rho_r$ restricted 
to the torsion part of $H^2(L^{p,q};\Z)$, i.e. 
\begin{eqnarray}\beta _r (\alpha)=u(\frac{p}{r}\Pi_{p,q}^* (x_r)+\frac{q}{r}\Pi_{p,q} ^* (y_r))\end{eqnarray}
for some $u\in (\Z /r)^*.$ 
We claim that modulo $r$ $b_1$ equals $u_r\frac{p}{r}$ and $b_2$ equals $u_r\frac{q}{r}.$ Thus $u=u_r$ for some (universal) $u_r\in (\Z/r)^*.$\\\\
\textbf{Proof of the last claim.} An idea to obtain information about the $\Pi_{p,q} ^* (x_r)$-component of $\beta _r (\alpha)$ is to analyze the \textit{restricted bundles}
 \[S^1\stackrel{i}{\rightarrow} L^{p,q}\vert _{S^2_1}\stackrel{\Pi_{p,q}^1}{\rightarrow}S^2\]
and
\[S^1\stackrel{j}{\rightarrow} L^{p,q}\vert _{S^2_2}\stackrel{\Pi_{p,q}^2}{\rightarrow}S^2,\]
where the first respectively the second fibre bundle is the restriction of the fibre bundle associated to $L^{p,q}$ to the first respectively the second $S^2$-factor. 
We realize that $L^{p,q}\vert _{S_1^1}$ and $L^{p,q}\vert _{S_2^2}$ are the familiar standard lens spaces $L_p^3$ and $L_q^3$ respectively.\\
Let $\alpha_p\in H^1(L_p^3;\Z /p )$ such that $\left\langle i^* (\alpha_p),[S^1]_{\Z /p}\right\rangle=1$ and $\alpha_q\in H^1(L_q^3;\Z /q )$ such that
$\left\langle j^* (\alpha_q),[U(1)]_{\Z /q}\right\rangle=1.$
Let $i_p$ and $i_q$ be the obvious inclusion of $L_p^3$ respectively $L_q^3$ in $L^{p,q}$ then it is clear that $i_p^*(\alpha)=:\alpha_{p,r}$ and $i_q^*(\alpha)=:\alpha_{q,r}$,
where $\alpha_{p,r}\in H^1(L^3_p;\Z/r)$ and $\alpha_{q,r}\in H^1(L^3_q;\Z/r)$ are the images of $\alpha_p$ and $\alpha_q$ respectively under the corresponding coefficient homomorphism.\\
Furthermore the following holds:
\begin{eqnarray*}&&i_p^*(\Pi_{p,q}^*(x))=\Pi^{1*}_{p,q}(x),  i_q^*(\Pi_{p,q}^*(x))=\Pi^{2*}_{p,q}(x), \\
 &&i_p^*(\beta_r(\alpha))=\beta_r (i_p^*(\alpha))=\beta_r(\alpha_{p,r}) \textrm{ and } i_q^*(\beta _r(\alpha))=\beta_r (i_q^*(\alpha))=\beta_r(\alpha_{q,r}).
\end{eqnarray*}
 We conclude from the construction of the maps together with the long exact sequence in $\Z/r$-cohomology
 for the pairs $(L^{p,q},L^3_p)$ and $(L^{p,q},L^3_q)$ that $i_p^*(\Pi_{p,q}^*(y_r))$ and $i_q^*(\Pi_{p,q}^*(x_r))$ vanish.\\
Summarizing the last considerations leads to the following:
\begin{eqnarray*}
&&i_p^*(\beta_r(\alpha))=i_p^*(u(\frac{p}{r}\Pi_{p,q}^* (x_r)+\frac{q}{r}\Pi_{p,q} ^* (y_r)))=u\frac{p}{r}\Pi^{1*}_{p,q} (x_r)=\beta_r(\alpha_{p,r}),\\
&&i_q^*(\beta_r(\alpha))=i_q^*(u(\frac{p}{r}\Pi_{p,q}^* (x_r)+\frac{q}{r}\Pi_{p,q} ^* (y_r)))=u\frac{p}{r}\Pi^{2*}_{p,q}(y_r)=\beta_r(\alpha_{q,r}).
\end{eqnarray*}
Thus if we knew $\beta_r(\alpha_{p,r})$ and  $\beta_r(\alpha_{q,r})$ in terms of $\Pi^{1*}_{p,q}(x_r)$ and $\Pi^{2*}_{p,q}(y_r)$ respectively then we would know what $\beta_r(\alpha)$ is.
\\
Assume $\beta_p(\alpha_p)=u_p\Pi_{p,q}^{1*}(x_p)$ for some $u_p\in(\Z/p)^*.$
We compare the short exact sequences associated to $\beta _r$ and $\beta _p$:\\
\setlength{\unitlength}{1cm}
\begin{picture}(10,3)
\put(2,2.45){$0$}
\put(2.4,2.55){\vector(1,0){0.9}}
\put(3.5,2.45){$\Z/r$}
\put(4.4,2.55){\vector(1,0){0.9}}
\put(3.9,1.6){\vector(0,1){0.5}}
\put(2.9,1.7){$\textrm{red}_{p,r}$}
\put(2,1.1){$0$}
\put(2.4,1.2){\vector(1,0){0.9}}
\put(3.5,1.1){$\Z/p$}
\put(4.4,1.2){\vector(1,0){0.9}}
\put(8.5,2.55){\vector(1,0){0.9}}
\put(8.5,1.2){\vector(1,0){0.9}}
\put(9.6,2.45){$0$}
\put(9.6,1.1){$0$,}
\put(5.5,2.45){$\Z/r^2$}
\put(6.5,2.55){\vector(1,0){0.9}}
\put(6.8,2.65){$\tilde{\pi}$}
\put(7.7,2.45){$\Z/r$}
\put(5.9,1.6){\vector(0,1){0.5}}
\put(4.8,1.7){$\textrm{red}_{p^2,r}$}
\put(8.1,1.6){\vector(0,1){0.5}}
\put(7.1,1.7){$\textrm{red}_{p,r}$}
\put(5.5,1.1){$\Z/p^2$}
\put(6.5,1.2){\vector(1,0){0.9}}
\put(6.8,1.3){$P$}
\put(7.7,1.1){$\Z/p$}
\end{picture}\\
where $\textrm{red}_{\cdot,\cdot}$ denotes the reduction homomorphism. The maps in the squares above commute, hence we get the following commutative diagram:
\\\\
\setlength{\unitlength}{1cm}
\begin{picture}(14,3)
\put(3.2,2.4){$H^1(L^3_p;\Z/r)$}
\put(5.7,2.5){\vector(1,0){1.5}}
\put(7.3,2.4){$H^2(L^3_p;\Z/r)$}
\put(6.3,2.6){$\beta_r$}
\put(3.2,0.5){$H^1(L^3_p;\Z/p)$}
\put(5.7,0.6){\vector(1,0){1.5}}
\put(7.3,0.5){$H^2(L^3_p;\Z/p)$,}
\put(6.3,0.75){$\beta_p$}
\put(4.3,1.0){\vector(0,1){1.1}}
\put(3.7,1.45){$\rho_{p,r}$}
\put(8.5,1.0){\vector(0,1){1.1}}
\put(8.6,1.45){$\frac{p}{r}\rho_{p,r}$}
\end{picture}\\\\
where $\rho_{\cdot,\cdot}$ is the "change of coefficient-homomorphism", i.e 
\[\beta _r \circ \rho_{p,r} = \frac{p}{r}\rho_{p,r} \circ \beta _p.\]
Thus $\beta_r(\alpha_{p,r})=u_{p,r}\frac{p}{r}\Pi^{1*}_{p,q}(x_r),$ where $u_{p,r}\in\Z/r$ is the mod-$r$-reduction of $u_p.$ It is clear that $u_{p,r}$ is a unit
 in $\Z/r.$ If $\beta_q(\alpha_{q,r})=u_q\Pi^{2*}_{p,q}(y_q)$ for some $u_q\in(\Z/q)^*$, then in the same way we obtain: 
$\beta_r(\alpha_{q,r})=u_{q,r}\frac{q}{r}\Pi^{2*}_{p,q}(y_r),$
where $u_{q,r}$ is the mod-$r$-reduction of $u_q.$ 
But (5) 
implies that $u_{p,r}=u_{q,r}=:u$ and thus $\beta _r (\alpha)=u(\frac {p}{r}\Pi_{p,q} ^* (x_r)+\frac{q}{r}\Pi_{p,q} ^* (y_r)).$
\\\\
By definition we have 
\begin{eqnarray*}f^*(v_1)&=&s\alpha,\\
f^*(\beta _r(v_1))&=&su_r(\frac {p}{r}\Pi_{p,q} ^* (x)+\frac{q}{r}\Pi_{p,q} ^* (y)),\\
f^*(z)&=&\epsilon(m\Pi_{p,q}^*(x)-n\Pi_{p,q} ^* (y))+k(\frac{p}{r}\Pi_{p,q}^*(x)+\frac{q}{r}\Pi_{p,q}^*(y)).
\end{eqnarray*}
Thus
\begin{eqnarray*}
f^*(v_1z^2)&=&-2s(\epsilon m+k\frac{p}{r})(\epsilon n-k\frac{q}{r})\alpha\Pi_{p,q}^*(xy),\\
f^*(v_1 (\beta _r v_1) z)&=&us((\epsilon m+k\frac{p}{r})\frac{q}{r}-(\epsilon n-k\frac{q}{r})\frac{p}{r})\alpha\Pi_{p,q}^*(xy),\\
f^*(v_1 (\beta _r v_1)^2 )&=&\frac{u^2s^3pq}{r^2}\alpha\Pi_{p,q}^*(xy).
\end{eqnarray*}
Since $\left\langle \alpha\Pi_{p,q}^*(xy),[L^{p,q}]_{\Z/r}\right\rangle$ equals $\left\langle \alpha,\Pi_{p,q}^*(xy)\cap [L^{p,q}]_{\Z/r}\right\rangle$ 
and by the choice of the orientation of $L^{p,q}$ it follows from \cite [Prop. 2] {G-77} and its proof that 
$\left\langle \alpha\Pi_{p,q}^*(xy),[L^{p,q}]_{\Z/r}\right\rangle\equiv 1 \textrm{mod }r$.
Hence since $2$ and $u_r$ are units of $\Z /r$ the equations (2)-(4) translate into the congruences stated in Theorem 8.\hfill $\square$\\\\
The manifolds in $\{L^{r,(t+lr)r}\}_l$ lie in the same simple and tangential homotopy type:\\\\
Here we have that $p=r$ and $q=(t+lr)r$ thus $\frac{p}{r}=1$ and $\frac{q}{r}=t+lr$. We choose $m$ to be 0 and $n$ to be 1. From Theorem 8 we obtain the following numbers modulo $r$:
\begin{eqnarray*}
sk(\epsilon-k(t+lr)),\textrm{ }
s^2((n+lr)k-(\epsilon-k(t+lr))),\textrm{ } s^3(t+lr).
\end{eqnarray*}
If we choose $s$ to be 1
and $k$ to be 0 then the three numbers above are modulo $r$ congruent to
\[0,\textrm{ }-\epsilon,\textrm{ }t\]
respectively which are independent of $l$. Thus the manifolds in $\{L^{r,(t+lr)r}\}_l$ 
lie in one homotopy type. By \cite [Lemma 12.5]{M-66} and Lemma 9 it follows that any homotopy equivalence between these manifolds is simple and
by the fact that any homotopy equivalence between closed manifolds is onto, Lemma 7 (ii) implies that these homotopy equivalences are also tangential.
Hence the manifolds in $\{L^{r,(t+lr)r}\}_l$ even lie in the same simple and tangential homotopy type.
\section{Distinguishing diffeomorphism types}
The diffeomorphism invariant which we use for smooth closed non-simply connected manifolds with finite cyclic fundamental group 
is the so called $\rho$-invariant which was introduced by M. Atiyah and I. Singer \cite{AS-68}. Let $M$ be a smooth closed non-simply connected and oriented
$5$-manifold with $\pi_1(M)\cong \Z/r$ being finite. Assume that there is a $6$-dimensional smooth oriented manifold $W$ with boundary $\tilde {M}$ which is equipped
with an orientation preserving smooth $\Z/r$-action such that the action coincides with the $\Z/r$-operation on the boundary $\tilde {M}$ given by deck transformation. 
Let $W_f$ be the fixed point set of the $\Z/r$-action. Assume the equivariant signature \cite{AS-68} of $W$ is trivial then the $\rho$-invariant of $M$ associated to
a non-trivial element $g$ of $\pi_1(M)$ is defined to be the evaluation of certain characteristic polynomials depending on the Chern-, Pontrjagin classes of 
the normal bundle of $W_f$ and the Pontrjagin classes of $W_f$, on the (twisted) fundamental class of $W_f$ if $W_f$ is orientable (not orientable).\\
Let $L^{p,q}\in \cal {L}$ and $r:=\vert\pi_1(L^{p,q})\vert$. Then we know that the universal covering space of $L^{p,q}$ is $L^{\frac{p}{r}\frac{q}{r}}$,
and that the deck transformation on $L^{\frac{p}{r}\frac{q}{r}}$ by $\pi_1(L^{p,q})$ is given by fibrewise rotation by angles corresponding to the $r$'th roots of unity. This 
perspective yields a canonical identification of $\pi_1(L^{p,q})$ with $\Z/r$. The disc bundle $D^{\frac{p}{r}\frac{q}{r}}$ associated to the $S^1$-fibre bundle structure
with the $\Z/r$-action canonically extended serves as a convinient choice of a bordism. Furthermore this $\Z/r$-bordism has trivial equivariant signature since
on the one hand the $\Z/r$-action is homotopically trivial, as it sits in an $S^1$-action, and on the other hand the dimension 
of the bordism is not divisible by 4 which means that the ordinary signature is trivial.\\
The fixed point set is just $S^2\times S^2$ and
the normal bundle of the fixed point set is isomorphic to the 2-dimensional real vector bundle given by the Euler class $\frac{p}{r}x+\frac{q}{r}y$.\\
Let $g$ be a non-trivial element of $\Z/r$
and $\theta_g$ the rotation angle between $0$ and $\pi$ of the action by $g$ then
\begin{eqnarray}\rho (g,L^{p,q})=-i\frac{\cos (\frac{\theta_g}{2})}{2r^2\sin^3 (\frac{\theta_g}{2})}pq,\end{eqnarray}
see e.g. \cite [p.88]{Ot-09}.
From this formula we see that the product of the parameters $p,q$ is a diffeomorphism invariant which shows that the manifolds in $\{L^{r,(n+lr)r}\}_l$ are
all pairwise non-diffeomorphic but by the previous section simply and tangentially homotopy equivalent which proves the first part of Theorem 1.
\section{Souls of codimension three}
Let $L^{p,q}\in \cal {L}$ and $(p,q,1)$ the matrix representation with respect to the standard basis of the epimorphism from $\Z^3$ to $\Z$ which sends $(x_1,x_2,x_3)$ to $px_1+qx_2+x_3$.
Let $\{\textbf{a}:=(a_1,a_2,a_3)$, $\textbf{b}:=(b_1,b_2,b_3)\}\in\Z^3$ be a basis of the kernel of this epimorphism
and let $i_{a,b}$ be the following Lie group homomorphism from $U(1)\times U(1)$ into $SU(2)\times SU(2)\times U(1)$:
\begin{eqnarray*}(z_1,z_2)&\mapsto& \left(\left(\begin{array}{cc}
z_1^{a_1}z_2^{b_1} & 0 \\ 
0&z_1^{-a_1}z_2^{-b_1}
\end{array} \right),\left(\begin{array}{cc}
z_1^{a_2}z_2^{b_2} & 0 \\ 
0 & z_1^{-a_2}z_2^{-b_2}
\end{array} \right),z_3^{a_3}z_2^{b_3}  \right).
\end{eqnarray*}
It is not hard to show that the image of $i_{a,b}$ does not depend on the choice of a basis of the kernel of $(p,q,1)$.
The following lemma proves the parts in the statements of Theorem 1 and 2 which involve the group $SU(2)\times SU(2)\times U(1)$ which we denote in this section by $G$.
\begin{lem}
The homogeneous space $\frac{G}{\mathrm{Im}(i_{a,b})}$ is diffeomorphic to $L^{p,q}$.
\end{lem}
\textbf{Proof.} The proof is subdivided into the following steps:
\begin{itemize}
 \item [(1)] First we define the projection map $\Pi:\frac{G}{\mathrm{Im}(i_{a,b})}\rightarrow S^2 \times S^2\times \C P^0$, where $\C P^0$ is a point
 and show that the fibre is $S^1$. Then we prove that there is a smooth and free $S^1$-action on $\frac{G}{\mathrm{Im}(i_{a,b})}$
 which preserves the fibre.
\item [(2)] We prove that there is a bundle isomorphism between $\Pi$ and the sphere bundle $S(E_{p,q})$ of the following complex line bundle:
\[E_{p,q}:(\textrm{pr}_1 ^* \gamma _2 ^p)\otimes (\textrm{pr}_2 ^* \gamma _1 ^q)\otimes \textrm{pr}_0^* \gamma _0 \rightarrow S^2 \times S^2\times \C P^0 ,\]
where $\gamma _2$ respectively $\gamma _1$ denotes the tautological bundle over the first and the second
factor of the base ($S^2=\C P^1$) respectively, $\gamma _0$ is the trivial complex line bundle $\C P^0$ and 
$\textrm{pr}_i:S^2 \rightarrow S^2 \times S^2\times \C P^0$ is the projection map onto the $i$-th factor ($i\in \{0,1,2\}$).
\end{itemize}
(1) Let us denote by $\dfrac{S^3 \times S^3 \times S^1}{\sim_{G_{a,b}}}$ the quotient of the following smooth $(S^1\times S^1)$-action on $S^3\times S^3\times S^1$:
\begin{eqnarray*}&&G_{a,b}:(S^1\times S^1)\times( S^3 \times S^3 \times S^1)\rightarrow S^3 \times S^3 \times S^1,\\
                 &&       (z_1,z_2,((x_1,x_2),(x_3,x_4),x_5))\mapsto (z_1^{a_1}z_2^{b_1}(x_1,x_2),z_1^{a_2}z_2^{b_2}(x_3,x_4),z_1^{a_3}z_2^{b_3}x_5)
\end{eqnarray*}
induced by the diffeomorphism $SU(2)\stackrel{\sim}{\rightarrow}S^3$, $A\mapsto A\cdot\left(\begin{array}{c} 1 \\ 0\end{array} \right)$. It is clear
that $\frac{S^3 \times S^3 \times S^1}{G_{a,b}}$ and $\frac{SU(2)\times SU(2)\times U(1)}{\mathrm{Im}(i_{a,b})}$ are diffeomorphic and we
denote this diffeomorphism by $\Phi$.
Let $\tilde{\Pi}$ be the map from $\dfrac{S^3 \times S^3 \times S^1}{\sim_{G_{a,b}}}$ to $S^2 \times S^2\times \C P^0$ 
which maps $[(x_1, x_2), (x_3,x_4),x_5]$ to $([x_1 :x_2], [x_3 : x_4],[x_5])$. We define $\Pi$ to be $\tilde{\Pi}\circ \Phi$.
\\\\
Let $d,e,f\in\Z$ such that $dp+eq+f=1$. There is the following split short exact sequence
\begin{eqnarray}1\rightarrow U(1)\times U(1)\stackrel{\left(\begin{array}{c}
a_1 \\
a_2 \\
a_3\end{array}\right. \left.\begin{array}{c} b_1\\b_2\\b_3  \end{array}  \right )_{\star}}{\longrightarrow}U(1)\times U(1)\times U(1) \begin{array}{c}\stackrel{(p,q,1)_{\star}}{\longrightarrow}\\
\stackrel{(d,e,f)_{\star}}{\longleftarrow}
\end{array}
 U(1)\rightarrow 1.
 \end{eqnarray}
From (7) we conclude that $\tilde{\Pi}$ is a $S^1$-fibre bundle, where the $S^1$-action is given by 

\[[(x_1,x_2),(x_3,x_4),x_5]\ast z =[(x_1,x_2)\cdot z^d, (x_3,x_4)\cdot z^e, x_5 \cdot z^f].\]
Hence $\Pi$ is a principal $U(1)$-fibre bundle.\\\\
(2) The following diagram commutes:

\setlength{\unitlength}{1cm}
\begin{picture}(10,3)
\put(2.3,2.42){$\dfrac{S^3 \times S^3 \times S^1}{\sim_{G_{a,b}}}$}
\put(7,2.4){$L^{p,q}$}
\put(4.5,0.5){$S^2\times S^2\times\C P^0\ \ ,$}

\put(6,2.6){$\psi$}
\put(5.7,2.5){\vector(1,0){1.1}}
\put(4.7,2){\vector(1,-1){1.06}}
\put(4.6,1.4){$\tilde{\Pi}$}
\put(7.4,2.3){\vector(-1,-1){1.355}}
\put(6.8,1.4){$\Pi_{p,q}$}

\end{picture}
\\
where \[\psi ([(x_1,x_2),(x_3,x_4),x_5])\]
equals
\[{\underbrace{\left( 
\begin{array}{cc}
	x_1\\x_2
\end{array}
\right)\otimes ... \otimes \left( 
\begin{array}{cc}
	x_1\\x_2
\end{array}
\right)}_{p \mbox{ copies}}\otimes\underbrace{\left( 
\begin{array}{cc}
	x_3\\x_4
\end{array}
\right)\otimes ... \otimes \left( 
\begin{array}{cc}
	x_3\\x_4
\end{array}
\right)}_{q \mbox{ copies}}}\otimes x_5\]
and
\begin{eqnarray*}
\Pi_{p,q} \left({\underbrace{\left( 
\begin{array}{cc}
	x_1\\x_2
\end{array}
\right)\otimes ... \otimes \left( 
\begin{array}{cc}
	x_1\\x_2
\end{array}
\right)}_{p \mbox{ copies}}\otimes\underbrace{\left( 
\begin{array}{cc}
	x_3\\x_4
\end{array}
\right)\otimes ... \otimes \left( 
\begin{array}{cc}
	x_3\\x_4
\end{array}
\right)}_{q \mbox{ copies}}}\otimes x_5\right ) \end{eqnarray*} equals \[([x_1:x_2], [x_3:x_4],[x_5]).\]
The assertion follows from the multiplicative property of the (total) Chern classes.\hfill $\square$\\\\
\textbf{Proof of Theorem 2.}
On the one hand the proof of \cite[Prop. 6.7.]{BKS-09} implies that tangentially homotopy equivalent manifolds,
 where the homotopy equivalences have trivial normal invariant, share the property that their cartesian product with $\R^3$ are diffeomorphic.
 On the other hand it shows that to any manifold $L^{r,qr}$ in $\{L^{r,(t+kr)r}\}_k$ there exists a subsequence of $\{L^{r,(t+kr)r}\}_k$ 
of manifolds which are homotopy equivalent, where the homotopy equivalences are tangential and have trivial normal invariants.
Thus by Theorem 1 there are infinitely many pairwise non-homeomorphic manifolds in $\cal {L}$ such that after taking the product with $\R^3$ they get diffeomorphic.
 The metrics we choose on each such 
product is the product metric of the submersion metric on the homogeneous quotient, where we choose the standard product metric on $SU (2)\times SU(2)\times U(1)$ and the 
euclidean metric on $\R^3$.\\
As an upper diameter bound we take the diameter $D$ of $SU(2)\times SU(2)\times U(1)$ with respect to the product metric of the standard metrics.
The existence of an uniform upper curvature bound follows from an idea of J.-H. Eschenburg \cite{E-72} and B. Totaro \cite {To-03} on $T^m$-actions which are 
subactions of isometric and free $T^{n}$-actions: Namely, O'Neill's formula shows that the
sectional curvature of a quotient manifold can be computed locally on $SU(2)\times SU(2)\times U(1)$. The same formula for the curvature formally makes sense for the
non-closed subgroup of $T^3\subset SU(2)\times SU(2)\times U(1)$ associated to any real linear subspace $\R^2$ in the Lie algebra $\R^3$ of $T^3$, where we use 
that $T^3$ acts freely on $SU(2)\times SU(2)\times U(1)$. The \textit{curvature} so defined is continuous on the compact manifold of all subspaces $\R^2$ of $\R^3$ and all 2-planes
in the tangent bundle of $SU(2)\times SU(2)\times U(1)$ which are orthogonal to the associated foliation of $SU(2)\times SU(2)\times U(1)$ by the $T^3$-action.
Hence there is a uniform upper bound for this curvature function and hence for the curvature of all quotients associated to subtori $T^2\subset T^3$.\hfill $\square$
\begin{rem}
An analysis of the surgery obstruction group, which lies beyond the scope of this work, shows that even between
any two manifolds in $\{L^{r,(t+kr)r}\}_k$ there exists simple and tangential homotopy equivalences with trivial normal invariant \cite [Corollary 19]{Ot-11}. Thus
all manifolds in $\{L^{r,(t+kr)r}\}_k$ get diffeomorphic after taking the product with $\R^3$.
\end{rem}
\begin{cor}
Let $r,t$ be as in Theorem 1 then the $\{L^{r,(t+kr)r}\}_k$ consists of pairwise (simply and tangentially) homotopy equivalent but non-diffeomorphic
Riemannian manifolds (the metrics are the homogeneous submersion metrics) with $0\leq sec\leq 1$ and diameter $\leq D$, where $D$ is a positive constant.
\end{cor}
This Corollary gives counterexamples in dimension 5 to the following mild version of a question by S.-T. Yau for nonnegative pinching \cite {Y-93}, Problem 11:\\\\
\textit{Let $n$ be a positive integer and $C,D>0$ real constants. Are there at most finitely many diffeomorphism classes of pairwise homotopy equivalent closed
Riemannian $n$-manifolds with $0\leq sec\leq C$ and diameter $\leq D$?}
\section{Remarks}
\begin{lem}
 There exists only one smooth transitive action of $SU(2)\times SU(2)\times U(1)$ on $L_r^3\times S^2$ up to self-diffeomorphisms of $L_r^3\times S^2$.
\end{lem}
\textbf{Proof.} Let $\Phi:SU(2)\times SU(2)\times U(1)\times L_r^3\times S^2\rightarrow L_r^3\times S^2$ be a smooth transitive Lie group action on $L_r^3\times S^2$. Then the 
isotropy group has to be 2-dimensional, connected and compact. But there is only $S^1\times S^1$ fulfilling these topological conditions. Thus up to Lie 
group automorphisms of $SU(2)\times SU(2)\times U(1)$ we are in the situation of the previous section, where we gave a canonical diffeomorphism between the 
homogeneous quotient and $L_r^3\times S^2$,
the total space of the principal $S^1$-fibre bundle over $S^2\times S^2$ given by the first Chern class $rx$ (see Introduction). But by formula (6) for the $\rho$-invariant
for spaces of such type we see that it is 0 which means that there are no other $S^1$-fibre bundles over $S^2\times S^2$ which have total spaces diffeomorphic to $L_r^3\times S^2$.
\hfill $\square$
\begin{lem}
\begin {itemize}
\item [(a)] Within the manifolds in $\cal {L}$ with odd order fundamental groups there doesn't exist pairs of non-homeomorphic manifolds which
can be realized as codimension 1 souls of a fixed manifold.
\item [(b)] Each sequence of pairwise distinct manifolds in $\cal {L}$ contains infinitely many $h$-cobordism classes.
\end {itemize}
\end{lem}
\textbf{Proof.} (a) Since the order of the fundamental groups of 
the manifolds under considerations is odd any real line bundle has to be trivial. Assume that there are two non-homeomorphic manifolds $L,L'$
which can be realized as codimension 1 souls of a manifold.
Then they are $h$-cobordant. Since the manifolds under consideration
have finite cyclic fundamental groups acting trivially on the cohomology of their universal covering space it 
follows from \cite {M-66} that the $h$-cobordism class is determined by the $R$-torisons of the boundary components. But from Lemma 9 it follows that
their Reidemeister torsions are trivial. The $s$-cobordism theorem 
implies that $L,L'$ have to be diffeomorphic which is a contradiction.\\\\
(b) This follows immediately from the calculation of the $\rho$-invariant in the previous section (formula 6) and the fact that this invariant is an $h$-cobordism invariant.
\hfill $\square$
\begin{rem}
\begin{itemize}
\item [(i)] From Lemma 19 and the proof of Theorem 2 it follows that there exists infinitely many pairwise non-equivalent smooth transitive actions of
 $SU(2)\times SU(2)\times U(1)\times \R^3$ on $L_r^3\times S^2\times \R^3$ with 
isotropy group $U(1)\times U(1)$, whereas there is only one class of smooth transitive $SU(2)\times SU(2)\times U(1)$-operations on $L_r^3\times S^2$.
 \item [(ii)] Lemma 20 (a) implies that there doesn't exist an infinite sequence of pairwise non-homeomorphic manifolds in $\cal {L}$ which can be realized as codimension 1 souls.
More general classes of lens space bundles over $S^2$ could deliver examples of pairs of non-homeomorphic manifolds which can be realized as codimension 1 submanifolds of an open manifold. Another problem is to prove the existence of metrics of nonnegative sectional curvature on such bundles.
\\
By Lemma 20 (b) there is not a sequence in $\cal {L}$ which consists of pairwise distinct manifolds such that only
finitely many $h$-cobordism classes occur. As a consequence of this observation we deduce from \cite[Prop. 6.10]{BKS-09} that there doesn't exist an infinite
sequence of distinct manifolds in $\cal {L}$ which can be realized as codimension 2 souls with trivial normal bundle of a fixed manifold.
\end{itemize}
\end{rem}

\end{document}